\documentclass[journal]{IEEEtran}

\usepackage{lineno,hyperref}
\modulolinenumbers[5]
\usepackage{amsmath,amsfonts,amssymb}
\usepackage{verbatim}
\usepackage{algpseudocode}
\usepackage{graphicx}
\usepackage{textcomp} %命令\textacutedbl的包,二阶导符号
\usepackage{multirow} %表格多行
\usepackage{epstopdf}
\usepackage{color}
\usepackage{algorithm}
\usepackage{threeparttable}
\usepackage[table,xcdraw]{xcolor}
\usepackage{subfigure}
\usepackage{changepage}
\usepackage{bm}
\usepackage{diagbox}

\newcommand{\At}{\mathbf{A}}
\newcommand{\St}{\mathbf{S}}

\newcommand{\Yt}{\mathbf{Y}}
\newcommand{\Nt}{\mathbf{N}}

\newcommand{\Pt}{\mathbf{P}}
\newcommand{\Qt}{\mathbf{Q}}
\newcommand{\Rt}{\mathbf{R}}

\newcommand{\phit}{\boldsymbol{\phi}}
\newcommand{\at}{\mathbf{a}}
\newcommand{\bbt}{\mathbf{b}}
\newcommand{\ct}{\mathbf{c}}
\newcommand{\dt}{\mathbf{d}}
\newcommand{\et}{\mathbf{e}}
\newcommand{\st}{\mathbf{s}}

\newcommand{\xt}{\mathbf{x}}

\newcommand{\mta}{\boldsymbol{\theta}}

\newcommand{\ft}{\mathbf{f}}

\newcommand{\ut}{\mathbf{u}}
\newcommand{\vt}{\mathbf{v}}

\begin{document}
\title{Gridless Evolutionary Approach for Line Spectral Estimation with Unknown Model Order}

\author{Bai~Yan,~
	Qi~Zhao,~	
	Jin~Zhang,~	
    J. Andrew~Zhang,~\IEEEmembership{Senior Member,~IEEE}~% <-this % stops a space
	and~Xin~Yao,~\IEEEmembership{Fellow,~IEEE}	
\thanks{Corresponding author: Jin Zhang.}
\thanks{B. Yan and Q. Zhao are with Guangdong Provincial Key Laboratory of Brain-Inspired Intelligent Computation, Department of Computer Science and Engineering, Southern University of Science and Technology, Shenzhen 518055, China, and also with School of Computer Science and Technology, University of Science and Technology of China, Hefei 230027, China (email: yanb@sustech.edu.cn; zhaoq@sustech.edu.cn).}% <-this % stops a space
\thanks{J. Zhang and X. Yao are with Guangdong Provincial Key Laboratory of Brain-Inspired Intelligent Computation, Department of Computer Science and Engineering, Southern University of Science and Technology, Shenzhen 518055, China (email: zhangj4@sustech.edu.cn; xiny@sustech.edu.cn).}
\thanks{J. A. Zhang is with Global Big Data Technologies Centre (GBDTC), University of Technology Sydney, NSW 2007, Australia (email: Andrew.Zhang@uts.edu.au).} }

% <-this % stops a space

\maketitle
\begin{abstract}
Gridless methods show great superiority in line spectral estimation. These methods need to solve an atomic $l_0$ norm (i.e., the continuous analog of $l_0$ norm) minimization problem to estimate frequencies and model order. Since this problem is \textit{NP-hard} to compute, relaxations of atomic $l_0$ norm, such as nuclear norm and reweighted atomic norm, have been employed for promoting sparsity. However, the relaxations give rise to a resolution limit, subsequently leading to biased model order and convergence error. 

To overcome the above shortcomings of relaxation, we propose a novel idea of simultaneously estimating the frequencies and model order by means of the atomic $l_0$ norm. To accomplish this idea, we build a multiobjective optimization model. The measurment error and the atomic $l_0$ norm are taken as the two optimization objectives. The proposed model directly exploits the model order via the atomic $l_0$ norm, thus breaking the resolution limit. We further design a variable-length evolutionary algorithm to solve the proposed model, which includes two innovations. One is a variable-length coding and search strategy. It flexibly codes and interactively searches diverse solutions with different model orders. These solutions act as steppingstones that help fully exploring the variable and open-ended frequency search space and provide extensive potentials towards the optima. Another innovation is a model order pruning mechanism, which heuristically prunes less contributive frequencies within the solutions, thus significantly enhancing convergence and diversity. Simulation results confirm the superiority of our approach in both frequency estimation and model order selection.
\end{abstract}

\IEEEpeerreviewmaketitle

\section{Introduction}
\IEEEPARstart{L}{ine} spectral estimation (LSE) aims at frequency estimation and model order selection from measurements collected as a superposition of complex sinusoids. Here, the ``model order selection'' means determining the number of frequencies. LSE has received significant attention as a major subject in signal processing fields. It has various applications, e.g., direction-of-arrival estimation in radar and sonar \cite{XiaJoint}, channel estimation in wireless communications, and simulation of atomic systems in molecular dynamics.  %参考2019年fast内点法for AST第2章

Many methods have been proposed for frequency estimation. Classical methods such as subspace methods \cite{RaoPerformance} are based on sample statistics. Their performance highly rely on a large number of snapshots and uncorrelated components. Moreover, the model order is required as a priori. With the development of compressive sensing theory, sparse methods have been presented for frequency estimation. These methods exhibit great advantages over sub-spaced methods, such as robustness to noise, no requirement for model order, and low requirement for snapshots. Spare methods can be divided into grid-based and gridless types. In the grid-based type \cite{2013A}\cite{YangOff}\cite{2013Spectral}\cite{2015Compressive}, the continuous frequency domain must be discretized into a finite grid, and the frequencies are restricted to this grid. Then the LSE problem is simplified to a sparse recovery problem. However, gridding gives rise to the well-known basis mismatch issue that limits the estimation accuracy. To avoid gridding, the so-called gridless type \cite{yang2014a}\cite{badiu2017variational}\cite{wagner2019gridless} have been proposed, which can directly operate in the continuous frequency domain. 

Apart from frequency estimation, the model order should also be determined in LSE. For subspace-based methods, some popular choices based on information criterion \cite{stoica2004model}\cite{valaee2004an}\cite{wax1985detection} can be incorporated to estimate the model order. It may be challenging to derive accurate model order in non-asymptotic regimes such as limited snapshots or low signal-to-noise-ratio (SNR) \cite{badiu2017variational}. Differently, sparse methods adopt sparse penalties to exploit frequencies' sparsity. Ideally, the $l_0$ norm or atomic $l_0$ norm is the best choice for sparse penalty, but they incur an \textit{NP-hard} problem. To make this problem easily solvable, the $l_0$ norm or atomic $l_0$ norm is usually relaxed to other sparse metric, e.g., $l_p$-norm ($p\in(0,1]$) sparse penalty, atomic norm \cite{Li2016Off}, reweighted atomic-norm \cite{yang2015enhancing}\cite{yang2018fast}, Gaussian prior \cite{Dai2017Sparse} or Gaussian-Bernoulli prior \cite{badiu2017variational}\cite{zhu2019grid}\cite{hansen2018superfast}. Such relaxation makes sparse methods suffer from a resolution limit \cite{yang2015enhancing}\cite{8963635}, i.e., the true frequencies are required to be well separated for successful recovery. Consequently, the resolution limit may lead to biased model order and frequency error.

In order to not only jointly estimate the model order and frequencies but also break the resolution limit, in this paper, we propose a multiobjective LSE model. The measurement error and atomic $l_0$ norm are taken as the two conflicting objectives. The multiobjective setting enables the frequencies and model order to be simultaneously estimated. Moreover, we exactly exploit the model order by the atomic $l_0$ norm without relaxation, thus breaking the resolution limit and providing accurate model order.  %we propose a multiobjective variable-length evolutionary search approach (MVESA) for gridless LSE. we propose a novel idea of simultaneously estimating the frequencies and model order by means of the atomic $l_0$ norm. To accomplish this idea,

The proposed multiobjective LSE model is a NP-hard problem. To solve this model, there are two issues to be addressed. One is how to handle the atomic $l_0$ norm without relaxation. Another is how to find the optima from the continuous (gridless) search space without knowing the true model order as a priori. 

Aiming at the two issues, we design a multiobjective variable-length evolutionary search algorithm (MVESA) to solve the proposed model. Here we design MVESA from the perspective of evolutionary algorithms, because evolutionary algorithms have revealed their strong ability to handle NP-hard $l_0$ problems \cite{yan2018adaptive}\cite{li2018preference}. To deal with the unknown model order problem, we introduce a variable-length coding and search strategy. This strategy flexibly codes diverse solutions with different lengths (i.e., different model orders). Then, it interactively searches diverse pathways (formed by solutions with different lengths) over the variable and open-ended frequency search space. These pathways act as steppingstones that help fully exploring the search space and provide extensive potentials towards the optima. Furthermore, we develop a model order pruning mechanism. This mechanism heuristically prunes less contributive frequencies within solutions. The pruning length is set at random. Hence, solutions' convergence and diversity is significantly improved. Finally, when the iterative generation terminates, the solution providing the most desirable trade-off between the two objectives is identified as the final solution.

Overall, this paper's main contributions are: 
\begin{itemize}
	\item Multiobjective LSE model. It simultaneously estimates the frequencies and model order without adjusting balancing parameters. Moreover, using atomic $l_0$ norm successfully breaks the resolution limit and provides accurate model order.
	\item Variable-length coding and search strategy. It flexibly codes solutions in different lengths and interactively searches diverse pathways over variable and open-ended search space, thus fast converging to the true frequencies. To our knowledge, it is the first time that realizes LSE over dynamic size of frequency search space. 
	\item Model order pruning mechanism. It heuristically prunes less contributive frequencies within solutions. The pruning length is set random. Therefore, solutions' convergence and diversity is both significantly improved.
	\item Empirical validation of MVESA's performance. Results confirm MVESA's efficacy and better performance in terms of frequency estimation and model order selection concerning to state-of-the-art methods. 
\end{itemize}

%\textbf{$\mathbb{R}$ and   }$\mathbb{E}$ denotes the mathematical expectation. $\At_{i,:}$ and $\At_{i,j}$ denote the $i$-th row and $(i,j)$-th element of a matrix $\At$.In particular, $\st|_\et$ stands for the sub-vector of $\st$ with entries indexed by the set $I=\{i|\et_i=1\}$. Similarly, $\St|_\et$ denotes the sub-matrix of $\St$ with rows indexed by the set $I=\{i|\et_i=1\}$. 

The rest of this paper is organized as follows. Section II provides background knowledge and related works. Sections III and IV present the proposed multiobjective atomic $l_0$ model and variable-length evolutionary search algorithm, respectively. Section V gives simulation results. Finally, section VI concludes the paper.

\textit{Notation:} Bold-face letters represent vectors and matrices, respectively. $\mathbb{R}$ and $\mathbb{C}$ denotes the real domain and complex one, respectively. $(\cdot)^T$, $(\cdot)^*$, and $(\cdot)^H$ denote transpose, conjugate, and conjugate transpose of a vector or matrix, respectively. 

\section{Background}
In this section, we introduce the LSE problem at first. Since our work falls into the gridless type and is closely related to evolutionary multiobjective optimization, we then provide a review of existing gridless methods and background knowledge of evolutionary multiobjective optimization.

\subsection{Line Spectral Estimation}
In LSE model, the measurements $\Yt\in\mathbb{C}^{M\times L}$ is a superposition of $K^*$ complex sinusoids corrupted by the white Gaussian noise $\Nt$:
\begin{equation}
	\begin{aligned}
		\Yt=\sum_{k=1}^{K^*}\at(\theta_k)\st_k^T+\Nt=\At\St+\Nt,
	\end{aligned}
	\label{eq-y=As+noise}
\end{equation}
where $\at(\theta_k)\triangleq[1,e^{j\pi\theta_k},...,e^{j(M-1)\pi\theta_k}]^T$ is the $k$-th complex sinusoid. $\theta_k\in[-1,1)$ and $\st_k\in\mathbb{C}^{L\times 1}$ denote the frequency and complex amplitudes of the $k$-th sinusoidal component. The $K^*$ complex sinusoids constitute $\At=[\at(\theta_1),\at(\theta_2),...,\at(\theta_{K^*})]\in\mathbb{C}^{M\times {K^*}}$. $\mta=[\theta_1,...,\theta_{K^*}]$ denotes the frequency combination. $\st_k^T$ is the $k$-th row of $\St$. The number of frequencies $K^*<M$, also referred as the ``model order'', is unknown in this paper. The goal of LSE is to estimate the model order
$K^*$ and frequencies $\mta$, given measurements $\Yt$ and mapping $\At_{(\mta)}$ (i.e., $\mta\rightarrow\At$).

\subsection{Related Works}
Gridless LSE methods do not need grid discretization and work directly in the continuous frequency domain. These methods need to solve an atomic $l_0$ norm (the continuous analog of $l_0$ norm) minimization problem. The atomic $l_0$ norm directly exploits sparsity and has no resolution limit, but it is NP-hard to compute. To make it tractable, earlier works switched to the convex atomic $l_1$ norm (also known as nuclear norm or atomic norm), including \cite{candes2014towards}\cite{tang2013compressed} for noiseless data and \cite{bhaskar2013atomic}\cite{yang2015gridless} for noisy data. Later, several works minimized a covariance matrix fitting criterion \cite{yang2014a}. They had been proved to be equivalent to atomic norm-based methods. However, due to the convex relaxation, the above methods suffer from a serious resolution limit, i.e., the frequencies are required to be well separated for recovery. 

To alleviate the resolution limit, the reweighted atomic-norm minimization \cite{yang2015enhancing}\cite{yang2018fast} and the reweighted covariance fitting criterion \cite{wu2018a} were reported to approximate the atomic $l_0$ norm. They brought enhanced sparsity and resolution compared to convex atomic norm-based methods. Alternatively, alternating projections-based gridless methods \cite{wagner2019gridless}\cite{wang2018ivdst:} directly solved the atomic $l_0$ norm minimization problem to pursue higher resolution. However, the convergence performance is not guaranteed due to unclosed or nonconvex sets\cite{wagner2019gridless}. Besides, by treating the frequencies as random variables, a few gridless LSE methods in Bayesian framework \cite{badiu2017variational}\cite{zhu2019grid}\cite{hansen2018superfast} were also proposed to estimate frequencies. 

Apart from frequency estimation, model order selection is also needed. Instead of using the atomic $l_0$ norm, atomic norm-based methods \cite{yang2015enhancing}\cite{yang2018fast}\cite{candes2014towards}\cite{yang2015gridless} exploit the model order by relaxed sparse metrics. This relaxation manner suffers from a resolution limit, subsequently producing biased model order and large frequency error. For covariance fitting criterion-based methods \cite{yang2014a}\cite{wu2018a}, the model order is usually identified by classic user-set threshold or information criterion methods \cite{stoica2004model}\cite{valaee2004an}\cite{wax1985detection} a posteriori. However, it is very challenging to derive accurate results due to inferior statistical properties in non-asymptotic regimes (e.g., limited snapshots or SNR) \cite{badiu2017variational}. For alternating projections-based methods, the model order is required as a priori \cite{wagner2019gridless}\cite{wang2018ivdst:}. In variable Bayesian methods \cite{badiu2017variational}\cite{zhu2019grid}\cite{hansen2018superfast}, the Gaussian-Bernoulli prior is employed to promote sparsity. However, it is still not yet clear how to determine the optimal sparse distributions in Bayesian framework \cite{8963635}.    

To summarize, with suboptimal/relaxed sparsity metrics or priors, existing gridless methods cannot faithfully promote sparsity and suffer from a resolution limit. Subsequently, the estimation accuracy of model order and frequencies is limited. Hence, it is expected to design a novel gridless method by means of the atomic $l_0$ norm without relaxation.

\subsection{Evolutionary Multiobjective Optimization}
Without loss of generality, we consider the multiobjective optimization problem (MOP) 
\begin{equation}
	\begin{aligned}
        &\min_{\xt} \ft(\xt)=(f_1(\xt),...,f_m(\xt)), \\&s.t.\ \xt\in\Omega 
	\end{aligned}
	\label{eq-MOPs}
\end{equation}
where $\xt$ is a candidate solution, $\Omega\subseteq\mathbb{R}^q$ is the search space, and $\ft$: $\Omega\rightarrow \mathbb{R}^m$ consists of $m$ real-valued objective functions. The objectives are conflicting to each other, which means no solution can minimize all the objectives simultaneously. 

\textbf{Definition 1.} Solution $\ut$ is said to Pareto dominate solution $\vt$, i.e., $\ut\prec\vt$, if and only if ${\forall}i\in\{1,2,...,m\}$, $f_i(\ut)\leqslant f_i(\vt)$, and ${\exists}j\in\{1,2,...,m\}$, $f_j(\ut)< f_j(\vt)$.   

\textbf{Definition 2.} $\xt^*$ is said to be a Pareto non-dominate solution, if there is no other solution $\xt\in\Omega$ satisfying $\xt\prec\xt^*$.

\textbf{Definition 3.} All the Pareto non-dominate solutions constitute Pareto optimal set, and their corresponding objective values form the Pareto front (PF). 

Evolutionary algorithms \cite{liu2019adaptive}\cite{Zhao2021Evolutionary} have been proposed to simultaneously optimize the multiple conflicting objectives in MOPs. These algorithms are available to various problem structures (e.g., non-convex, non-linear) and can obtain a set of Pareto non-dominate solutions with different trade-off among the multiple objectives in a single run.   

\section{Proposed Multiobjective LSE Model}
To simultaneously estimate frequencies and model order, we naturally formulate the LSE model (\ref{eq-y=As+noise}) as an MOP. The measurment error and the atomic $l_0$ norm are taken as two conflicting objectives. Our formulation holds two advantages: 1) frequencies and model order can be simultaneously estimated without adjusting any balancing parameter; 2) the model order is exactly exploited by the atomic $l_0$ norm without relaxations, hence breaking the resolution limit. 

For clarity, we first profile the atomic $l_0$ norm before giving our proposed model. We follow the research \cite{yang2016exact} to define the atomic $l_0$ norm of measurements $\Yt$. Specifically, define an atomic set
\begin{equation}
\begin{aligned}
\mathcal{A}:=\{\dot{\at}(\theta,\phit):=\at(\theta)\phit: \theta\in[-1,1), \phit\in\mathbb{C}^{1\times L},\|\phit\|_2=1\},
\end{aligned}
\label{eq-atomic set}
\end{equation}
it can be viewed as an infinite dictionary indexed by the continuous varying parameters $\theta$ and $\phit$. The atomic $l_0$ norm of measurements $\Yt$, $\|\Yt\|_{\mathcal{A},0}$, is defined as the minimum number of atoms in $\mathcal{A}$ that synthesizes $\Yt$: 
\begin{equation}
\begin{aligned}
\|\Yt\|_{\mathcal{A},0}=\inf_{\theta_k, \phit_k, c_k} &\{\kappa: \Yt=\sum_{k=1}^{\kappa}\dot{\at}(\theta_k, \phit_k)c_k, \theta\in[-1,1), \\&\qquad\|\phit\|_2=1, c_k>0\},\\=\inf_{\theta_k, \st_k} &\{\kappa: \Yt=\sum_{k=1}^{\kappa}\at(\theta_k)\st_k, \theta\in[-1,1)\},
\end{aligned}
\label{eq-atomic l0 norm}
\end{equation}
where ``inf" stands for infimum, $\phit_k=c_k^{-1}\st_k$, and $c_k=\|\st_k\|_2$.

By introducing the atomic $l_0$ norm (\ref{eq-atomic l0 norm}), we formulate LSE as a multiobjective optimization model
\begin{equation}
\begin{aligned}
&\min \ft(\mta, \St)=(\|\Yt\|_{\mathcal{A},0}, \|\Yt-\At\St\|_F^2),
\end{aligned}
\label{eq-MOP}
\end{equation}
where the two objectives, $\|\Yt\|_{\mathcal{A},0}$ and $\|\Yt-\At\St\|_F^2$, represent the atomic $l_0$ norm and measurement error, respectively. Each solution is composed of a frequency combination $\mta$ and amplitudes $\St$. A major advantage of this model is that, the exact atomic $l_0$ norm is introduced to appropriately promote sparsity, thus it does not suffer from a resolution limit compared to existing gridless methods \cite{wu2018a}\cite{wang2018ivdst:}. Hence, estimates of model order and frequencies can be more accurate.% It is evident that the proposed model (\ref{eq-MOP}) is NP-hard due to the $l_0$ like norm. 

\section{Proposed Variable-length Evolutionary Search Algorithm}
We design a variable-length evolutionary search algorithm to solve the proposed multiobjective LSE model (\ref{eq-MOP}). The designed algorithm includes two major innovations. One is a variable-length coding and search strategy. It flexibly codes and interactively searches diverse solutions with different model orders. These solutions act as steppingstones that help fully exploring the variable and open-ended frequency search space and provide extensive potentials towards the optima. Another innovation is a model order pruning mechanism. It heuristically prunes less contributive frequencies within the solutions. The pruning length is set at random. With this mechanism, solutions' convergence and diversity can be significantly enhanced.

\begin{algorithm}[t]
	\caption{Pseudo-code of MVESA} 
	\begin{algorithmic}[1]
		\Require mapping $\At_{(\mta)}$, measurements $\Yt$, empty archive $\Rt^G$
		\Ensure	$\mta$
		\State $G=1$;
		\State $\Pt^G\leftarrow$\emph{Initialization}$(\At_{(\mta)},\Yt)$;	
		\While {``\textit{stopping criterion not met}''}		
		\State $\Pt^G\leftarrow$\emph{Tournament\_Selection}$(\Pt^G)$;			
		\State $\Qt\leftarrow$\emph{Variable-length\_Search}$(\Pt^G)$;			
		\State $\Pt^G\leftarrow$\emph{Environmental\_Selection}$(\Pt^G\bigcup\Qt, \At_{(\mta)},\Yt)$;
		\State $(\Pt^G, \Rt^G)\leftarrow$\emph{Archiving\&Model\_Order\_Pruning}$(\Pt^G,$ $\Rt^G, \At_{(\mta)}, \Yt)$; 		
		\State $G=G+1$;	
		\EndWhile	
		\State $\mta\leftarrow$\emph{Knee\_Solution\_Identification}$(\Rt^G)$;			
	\end{algorithmic}
	\label{al-MVESA}
\end{algorithm}

\subsection{Overall Framework}
The workflow of the proposed MVESA is shown in Algorithm \ref{al-MVESA}. MVESA starts with initialization. A number of $N$ frequency combinations with different model orders are generated at random. Their corresponding amplitudes are recovered via the least square method (will be detailed in equation (\ref{eq-decoding}), Section \ref{sec-coding}). The $N$ frequency combinations and amplitudes compose the initial solution set $\Pt^G$, where $G$ is the generation counter. $\Pt^G$ is evaluated by model (\ref{eq-MOP}). Iterative generations follow the initialization. In each generation, the following steps are executed:

\textbf{Step 1}, tournament selection. The well-known binary tournament selection operator \cite{blickle2000tournament} is employed to select more effective initial solutions as parents. These parents will be used to produce offspring solutions in the next step.   %does not guarantee reproduction of best solution

\textbf{Step 2}, variable-length search. With obtained parents, a modified variable-length crossover (will be detailed in Section \ref{sec-coding}) and the polynomial mutation operator \cite{deb1996combined} are introduced to produce $N$ new frequency combinations with different model orders. New frequency combinations' amplitudes are recovered by the least square method. The new frequency combinations and their amplitudes make up the offspring solution set $\Qt^G$. $\Qt^G$'s fitness is calculated by model (\ref{eq-MOP}).

\textbf{Step 3}, environmental selection. The environmental selection operator of NSGA-II \cite{deb2002fast} is applied to select $N$ better (in terms of convergence and diversity) solutions from $\Pt^G\bigcup\Qt$. These $N$ solutions reform $\Pt^G$.

\textbf{Step 4}, archiving and model order pruning. We set an external archive $\Rt^G$ to collect the best solutions with each possible model order found so far (i.e., from $\bigcup_{G=1}^G\Pt^G$). This archive avoids missing optimal solutions during iterations. $\Rt^G$ is firstly updated with elite solutions $\Pt^G$. The solutions that newly join $\Rt^G$ at the current generation are denoted as \textit{newcomers}. We design a novel model order pruning mechanism to significantly improve the archive's convergence and diversity. Specifically, for each newcomer, this mechanism heuristically prunes its less contributive frequencies. Then, the pruned newcomer's amplitudes are recovered by the least square method. With pruned newcomers' frequency combiantions and amplitudes, the fitness of pruned newcomers is obtained by model (\ref{eq-MOP}). Finally, archive $\Rt^G$ and population $\Pt^G$ are updated with the pruned newcomers. 

Once the iterative generation terminates, we identify the knee solution from archive $\Rt$ as the final solution. This is because the knee solution has the maximum marginal rates of return, that is, an improvement in one objective would lead to a severe degradation in another. It provides an attractive trade-off between the two objectives \cite{RachmawatiMultiobjective}, and the efficacy has been empirically confirmed in Section \ref{sec-exp-knee}. Here we employ the kink method \cite{mierswa2006information} to identify the knee solution, by which the solution with the largest slope variance over the obtained PF is taken as the knee. Finally, MVESA returns the knee solution's frequency combination as the output. 

Core components of MVESA, i.e., the variable-length coding and search, archiving and model order pruning, are detailed below.

\subsection{Variable-length Coding and Search}\label{sec-coding}
The true model order is unknown in advance. Obtrusively using a predefined fixed-length (fixed model order) coding would lead to sub-optimal, deteriorating LSE performance. To handle this, we naturally introduce a variable-length coding strategy to represent solutions with diverse possible model orders, offering tremendous flexibility. To our knowledge, it is the first time that realizes direct LSE over dynamic size of frequency search space.  %such that we can fully explore  the variable, open-ended search space to find the optimal solution. 

\begin{figure}[t]
	\centering
	\includegraphics[width=5.5cm,height=2.2cm]{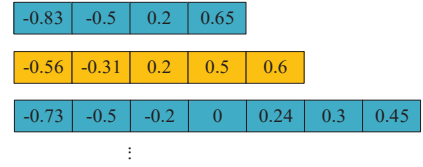}
	\caption{Variable-length coding of frequency combinations. Each row represents a frequency combination.}
	\label{fig-coding}
\end{figure}

In the variable-length coding strategy, we code each solution by a variable-length representation 
\begin{equation}
	\begin{aligned}
		\Pt&=\{(\mta_1, \St_1),...(\mta_n, \St_n),...,(\mta_N, \St_N)\}, \\ \mta_n&=[\theta_{n1},\theta_{n2},...,\theta_{nd_n}]\in\mathbb{R}^{1\times d_n}, \\
		\St_n&\in\mathbb{C}^{d_n\times L},
	\end{aligned}
	\label{eq-val-representation}
\end{equation} 
where the population $\Pt$ consists of $N$ solutions, $(\mta_n, \St_n)$ is the $n$-th solution, and $d_n$ is the length (model oder) of the $n$-th solution. The variable-length is reflected by solutions with different $d_n$s. The model order could be evolved towards the ground-truth during evolutionary search. Fig. \ref{fig-coding} gives an example of frequency combinations. Each row indicates a frequency combination, where frequencies are sorted in an ascending order, and the number of lattices is the model order.  %the decimal in each lattice stands for a frequency;

Once a $\mta$ is obtained by the initialization or evolutionary search, we can employ a simple-yet-effective least square method to solve $\min_{\St} \|\Yt-\At\St\|_2$ and acquire the corresponding amplitudes $\St$:  
\begin{equation}
\begin{aligned}
\St=(\At^T\At)^{-1}\At^T\Yt.
\end{aligned}
\label{eq-decoding}
\end{equation} 
Therefore, the task of LSE becomes to find the frequency combination as accurately as possible. 
%Crossover and mutation play a crucial rule in producing offspring. Therefore, we need to carefully set suitable crossover and mutation operators based on the variable-length representation. 

To find the optimal frequency combination, we should discreetly design evolutionary search operators. Since mutation have no concern with solution's lengths, we employ the classical polynomial mutation \cite{deb1996combined} to perturb solutions. However, traditional crossover operators cannot be incorporated because they are only designed for fixed-length coding. Recently, quite a few variable-length crossover operators have been proposed \cite{ryerkerk2019survey}, e.g., cut and splice, spatial, and similarity-based operators. Cut and splice operators are the most disruptive. Spatial operators are the least disruptive but can only be applied to specific problems with spatial components. Similarity-based operators are less disruptive by preserving common sequences and allowing only differences to be exchanged or removed. Thus, we modify a similarity-based operator, i.e., the synapsing variable-length crossover \cite{hutt2007synapsing}, and incorporate it into our work. 

\begin{figure}[t]
	\centering
	\includegraphics[width=8cm,height=8cm]{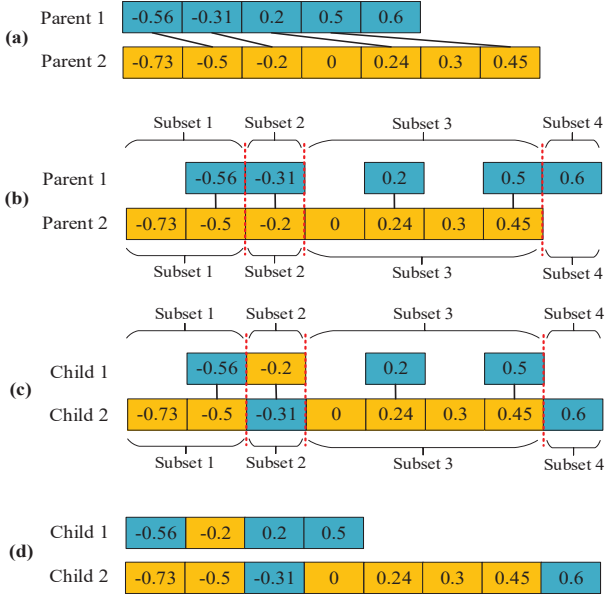}
	\caption{Variable-length crossover. (a) Link the most similar counterparts between the two parents by black oblique lines. (b) Align the two parents, and segment the two parents into $\bar{n}+1$ paired subsets by $\bar{n}$ red dotted lines, where $\bar{n}=3$ is the number of crossover points. (c) Crossover within each even paired subsets, respectively. (d) Glue subsets to produce two offsprings.}
	\label{fig-crossover}
\end{figure}

%We aim to maintain solutions with diversified model orders, such that the solutions can act as steppingstones helping fully exploring the search space and providing extensive potentials towards the true frequencies.
We modify the synapsing variable-length crossover by randomly selecting the number of crossover points. This modification helps maintain good diversity of solutions. The process is exhibited in Fig. \ref{fig-crossover}, which involves two steps:

1) Aligning parents. As displayed in Fig. \ref{fig-crossover}(a), black oblique lines link the lattice of one parent to a counterpart of another parent, such that the linked two lattices are the most similar with respect to each other. The similarity is measured by the Euclidean distance. Based on the links, we can align the two parents, as shown in Fig. \ref{fig-crossover}(b).   

2) Executing $\bar{n}$-point crossover. We randomly choose the number of crossover points $\bar{n}$ for maintaining diversity. Here, $\bar{n}$ is a random integer value between 1 and the length of the shorter parent. With $\bar{n}$ crossover points (red dotted lines in Fig. \ref{fig-crossover}(b)), the two parents are segmented into $\bar{n}+1$ paired subsets. According to the principle of $\bar{n}$-point crossover, each even paired subsets exchange affiliations, respectively. For example, in Fig. \ref{fig-crossover}(c), the second paired subsets $\{-0.31\}$ and $\{-0.2\}$ exchange their affiliations; so does the fourth paired subsets $\{0.6\}$ and $\{ \}$. Finally, subsets are glued to produce two offsprings with different lengths, as depicted in Fig. \ref{fig-crossover}(d). 

\subsection{Archiving and Model Order Pruning Mechanism} \label{sec-archive_pruning}
We propose a novel archiving and model order pruning mechanism, which tremendously enhances solutions' convergence and diversity. Archiving refers to using the external archive $\Rt$ to collect the best solutions with each possible model order found during iteration. It avoids missing optimal solutions. Furthermore, the convergence and diversity of archive solutions can be well maintained. Model order pruning aims to prune less contributive frequencies within solutions, thus reducing solutions' redundancy andredundancy enhancing the convergence perofrmance.

The motivation of designing model order pruning step is as follows. The variable-length search is very likely to produce long solutions. These solutions may include both close-to-optimal frequencies and spurious ones. It is necessary to prune the spurious frequencies and push overlong solutions towards the optima. Generally, close-to-optimal frequencies possess higher power than the spurious ones do. Motivated by this, we design the model order pruning mechanism to heuristically get rid of frequencies with lower power, so that the frequency combination's redundancy is greatly reduced and the resulted solution length approaches the true model order.

The pseudo-code of archiving and model order pruning is shown in Algorithm \ref{al-archiving and refinement}. It includes three operations: archiving, model order pruning, and update. 

\begin{figure*}[t] 
	\centering
	\subfigure[Archiving.]{\includegraphics[width=0.44\textwidth]{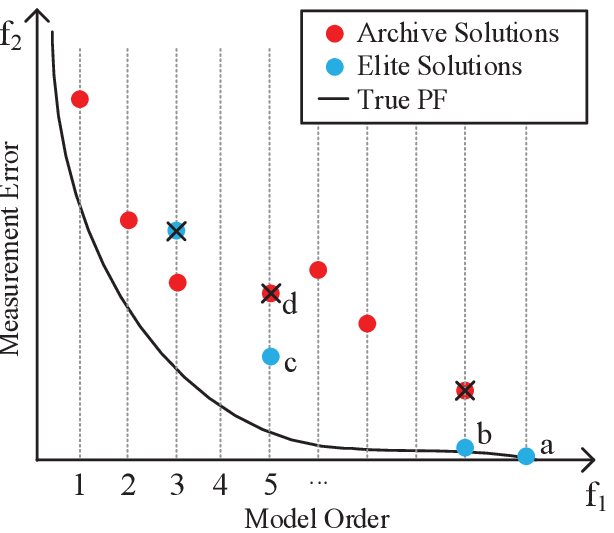}} \qquad
	\subfigure[Model order pruning.]{\includegraphics[width=0.44\textwidth]{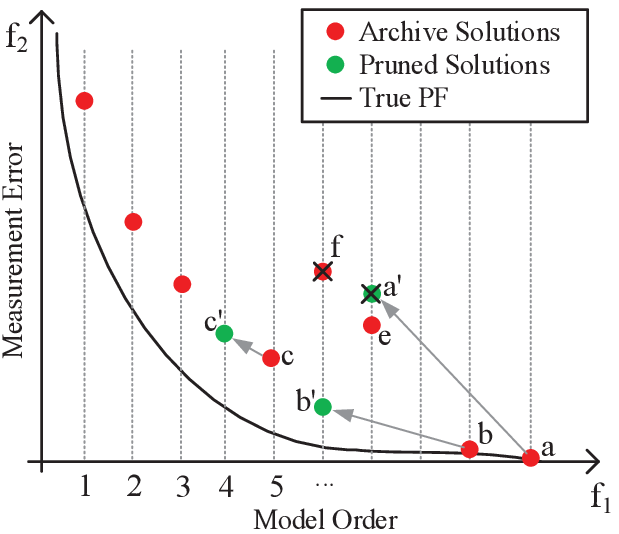}}
	\caption{Archiving and model order pruning mechanism. Vertical dotted lines refer to different model orders. Points labeled with cross symbols indicate being discarded.}
	\label{fig-refinement}
\end{figure*}

\textbf{Archiving (lines 2-5 of Algorithm \ref{al-archiving and refinement})}. We aim to collect the best solutions with each possible model order so far and store them into archive $\Rt^G$. For clarity, we denote the Pareto non-dominate solutions of population $\Pt^G$ as \textit{elite solutions}. As depicted in Fig. \ref{fig-refinement}(a), we put the archive solutions and elite solutions together. For each pair of archive and elite solutions with the same length, we replace the archive solution by the elite one only if this archive solution is dominated by the elite one. Thereafter, the solutions that newly join the archive are denoted as \textit{newcomers}. %

For example, in Fig. \ref{fig-refinement}(a), the \textit{archive} solution $\dt$ and the \textit{elite} solution $\ct$ has the same model order. $\dt$ is dominated by $\ct$, so we replace $\dt$ in archive by $\ct$, enhancing the convergence performance. In this way, we determine all newcomers $\{\at,\bbt,\ct\}$. With the archiving mechanism, solutions with the best performance for each possible length can always be reserved, thus the convergence and diversity of solutions is both improved.

\begin{algorithm}[t]
	\caption{\emph{Archiving\&Model\_Order\_Pruning}}
	\begin{algorithmic}[1]	
		\Require population $\Pt^G$, archive $\Rt^G$, mapping $\At_{(\mta)}$, measurements $\Yt$
		\Ensure updated population $\Pt^G$, updated archive $\Rt^G$
		\State /*Archiving*/ 		
		\For{``each possible length $K$"}		
		\State Replace the archive solution with length $k$ by the elite solution with the same length, only if this archive solution is dominated by this elite solution;		
		\EndFor
		\State Identify the \textit{newcomers} in $\Rt^G$;
		\State /*Model order pruning*/
		\For{each \textit{newcomer}}
		\State Compute each frequency's power via equation (\ref{eq-frequency power});
		\State Sort the frequencies' power in descending order;
		\State Acquire the pruned frequency combination based on equation (\ref{eq-prune frequency}).
		\EndFor
		\State Obtain the pruned \textit{newcomers} $(\mta_{Trun},\St_{Trun})$ via equation (\ref{eq-decoding}) and compute its fitness; 			
		\State /*Update*/			
		\State Update population $\Pt^G$ and archive $\Rt^G$ with \textit{newcomers} according to the three updating cases. 
	\end{algorithmic}
	\label{al-archiving and refinement}
\end{algorithm}

\textbf{Model order pruning (lines 7-12 of Algorithm \ref{al-archiving and refinement})}. We execute this operation to reduce the redundancy of newcomers, as exhibited in Fig. \ref{fig-refinement}(b). Assume $\xt\in\{\at,\bbt,\ct\}$ is a newcomer, $\xt$ corresponds to the decision variable $(\mta,\St)$, and the model order of $\xt$ is $\hat{K}$, we prune the frequency combination $\mta$ by three steps:

1) Calculating each frequency's power by averaging the power over multiple snapshots 
\begin{equation}
\begin{aligned}
p_i=\sqrt{\sum_{l=1}^{L}|S_{il}|^2}, i=1,2,...,\hat{K},
\end{aligned}
\label{eq-frequency power}
\end{equation}
where $p_i$ is the $i$-th component's power of newcomer $\xt$, $ S_{il}$ is the $i$-th row and $l$-th column element of $\St$.

2) Sorting the $\hat{K}$ frequencies' powers in descending order: 
\begin{equation}
\begin{aligned}
p_{i_1}\geqslant p_{i_2}\geqslant...\geqslant p_{i_{\hat{K}}},
\end{aligned}
\label{eq-sort frequency power}
\end{equation}
where $i_1$,$i_2$,...,$i_{\hat{K}}$ is a permutation of $\{1,2,...,\hat{K}\}$.

3) Pruning the frequency combination $\mta$ and retaining the frequencies with high power in priority. To maintain the diversity in length, the length to be cut off $\hat{K}_{cut}$ is set to be a random value from $[1,\hat{K}-1]$. After cutting off, $\mta$ becomes
\begin{equation}
	\begin{aligned}
		\mta=[\theta_{i_1},\theta_{i_2},...,\theta_{i_{\hat{K}-\hat{K}_{cut}}}],
	\end{aligned}
\label{eq-prune frequency}
\end{equation}

As observed in Fig. \ref{fig-refinement}(b), with model order pruning, the newcomers $\{\at,\bbt,\ct\}$ give birth to pruned solutions $\{\at',\bbt',\ct'\}$. For each pruned solution, the corresponding amplitudes is recovered by equation (\ref{eq-decoding}), and the fittness can be obtained via model (\ref{eq-MOP}). 

\textbf{Update (lines 14 of Algorithm \ref{al-archiving and refinement})}. With pruned solutions, we update the archive $\Rt^G$ and population $\Pt^G$, as shown in Fig. \ref{fig-refinement}(b). One of the following three updating cases would occur: %Points $\{a',b',c'\}$ are the produced truncated solutions based on newcomers $\{a,b,c\}$, respectively.
\begin{itemize}
\item Case $\at'$: The pruned solution $\at'$ is dominated by the archive solution $\et$ with the same length, i.e., $\et\prec\at'$. Therefore, $\Rt^G$ and $\Pt^G$ remain unchanged. 

\item Case $\bbt'$: The pruned solution $\bbt'$ dominates the archive solution $\ft$ with the same length, i.e., $\bbt'\prec\ft$. Therefore, $\bbt'$ substitutes $\ft$ in $\Rt^G$ and substitutes an arbitrary solution in $\Pt^G$. 

\item Case $\ct'$: The pruned solution $\ct'$ is generated via pruning the newcomer $\ct$. There exists no archive solution with the same length as $\ct'$. Therefore, $\ct'$ is added to $\Rt^G$.
\end{itemize}
It can be observed that, with archiving and model order pruning mechanism, the resulted archive (i.e., points without a cross in Fig. \ref{fig-refinement}(b)) can obtain better convergence and diversity performance.

\subsection{Convergence and Complexity Analysis}\label{sec-complexity}
We now analyze the convergence and computational complexity of the proposed MVESA.

\textit{1) Convergence:} The convergence of proposed MVESA can be guaranteed.

\textit{Analysis:} In MVESA, the archive $\Rt$ can be viewed as reserving the best solutions to a series of subproblems with different model orders. Specifically, for a subproblem associated with a specific model order, $\Rt$ always reserves a solution with lower measurement error to this subproblem. It implies that for each subproblem, the measurement error is non-increasing after each iteration of MVESA. Since the measurement error is lower bounded for each subproblem, MVESA is guaranteed to converge. $\hfill\blacksquare$

\textit{2) Complexity:} The main computational complexity of MVESA lies in the modified crossover and amplitudes recovery. At each generation, the modified crossover requires $\mathcal{O}(NM^2)$ computations, where $N$ and $M$ are the population size and the number of measurements, respectively. The computational complexity of recovering amplitudes is $\mathcal{O}(2NM^3)$ in the worst case. Thus, the total complexity is $\mathcal{O}(NM^3)$.   

\section{Simulation Experiments}\label{sec-simulation}
In this section, we conduct simulation experiments to evaluate the performance of the proposed MVESA compared to state-of-the-art algorithms under various scenarios.

%下面两个图属于下一个subsection
\begin{figure*}[t] 
	\centering
	\subfigure[]{\includegraphics[width=0.46\textwidth]{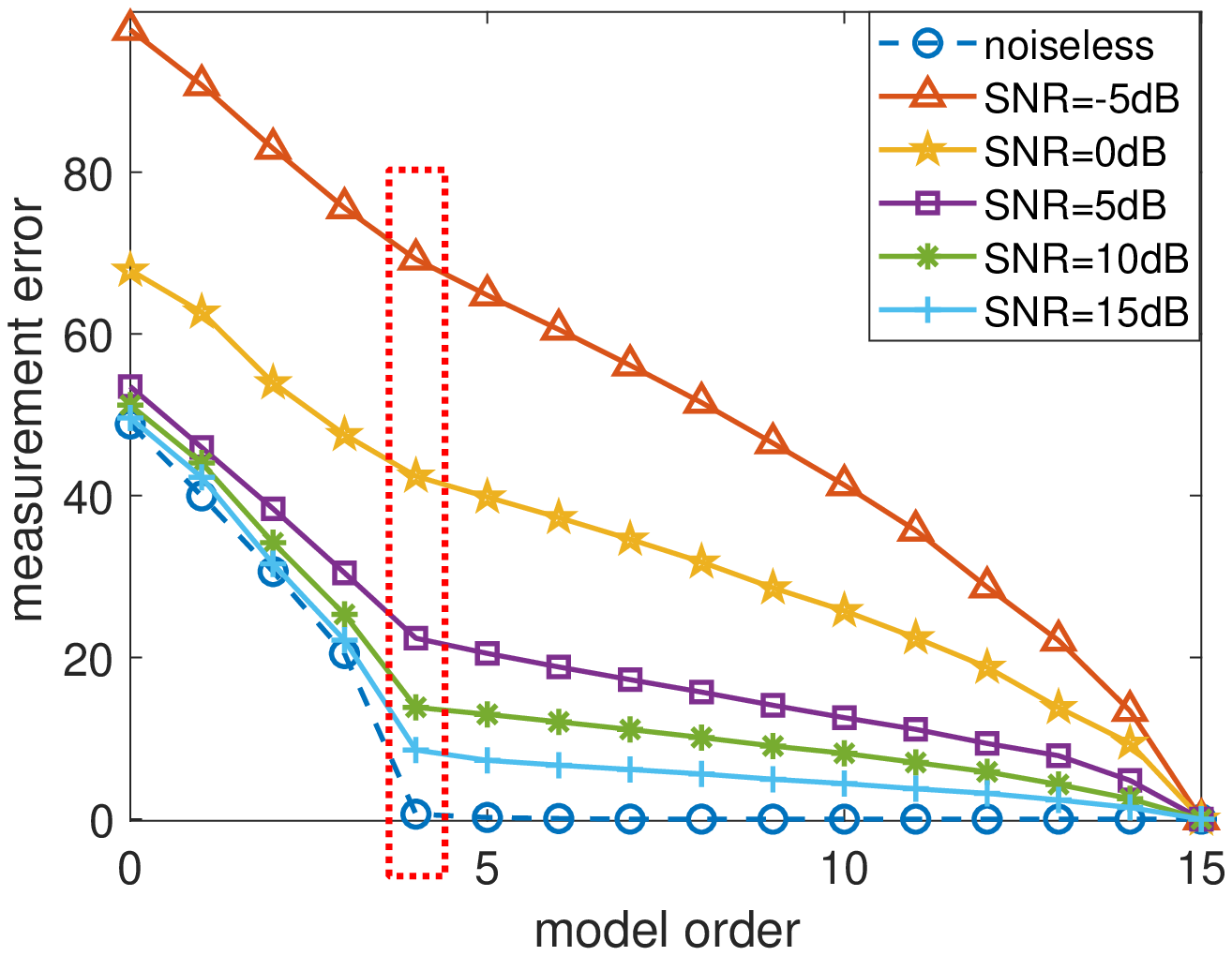}} \qquad
	\subfigure[]{\includegraphics[width=0.43\textwidth]{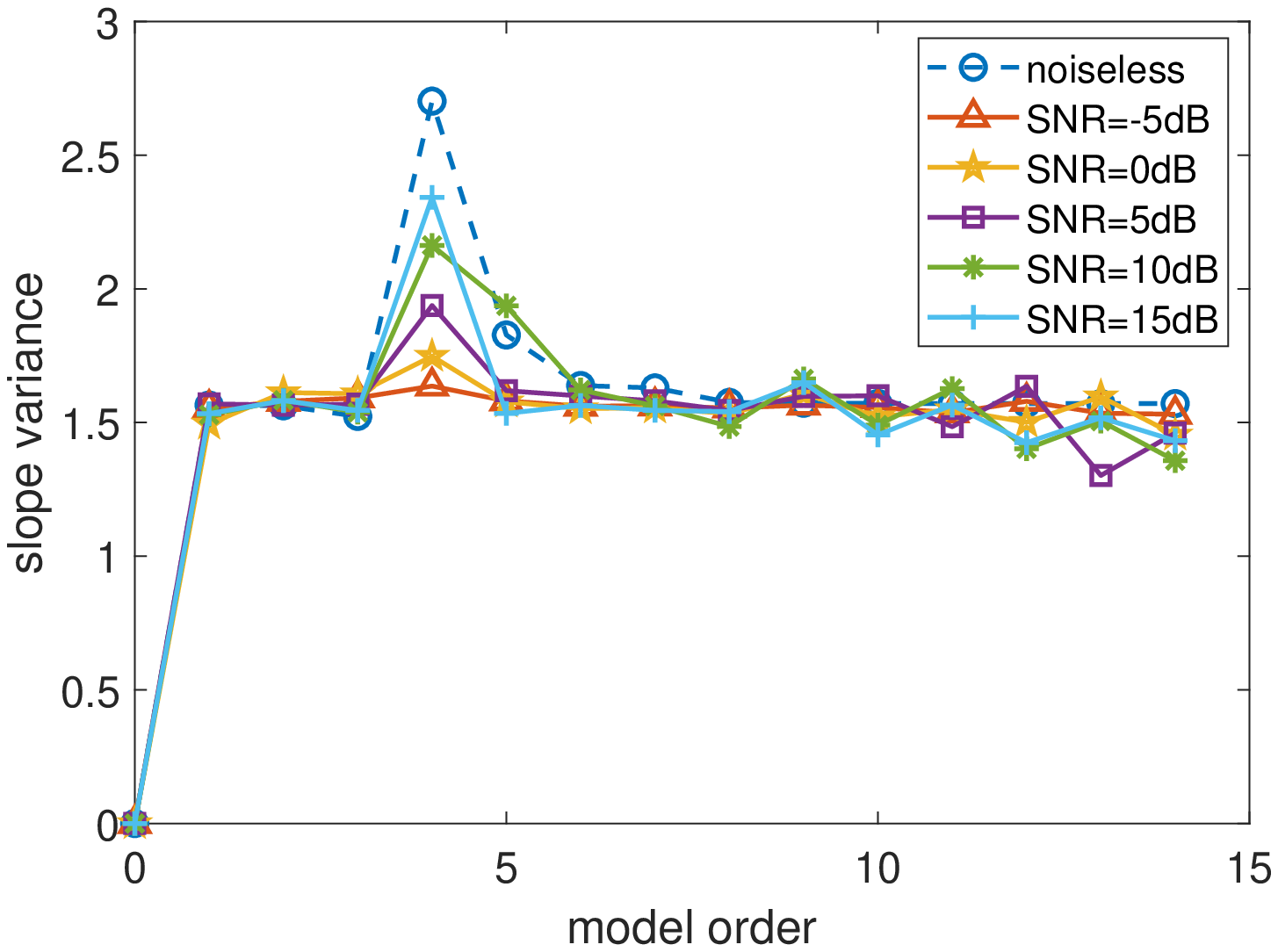}}
	\caption{Pareto front (a) and slope variance (b) of final archive versus SNR for $K=4$, $M=15$, and $T=20$. The points inside the dotted line denote identified knee solutions.} %theta=[-16.7  6.4 14.9 24.5];
	\label{fig-PF-and-slope}
\end{figure*} 
\begin{figure*}[t] 
	\centering
	\subfigure[]{\includegraphics[width=0.42\textwidth]{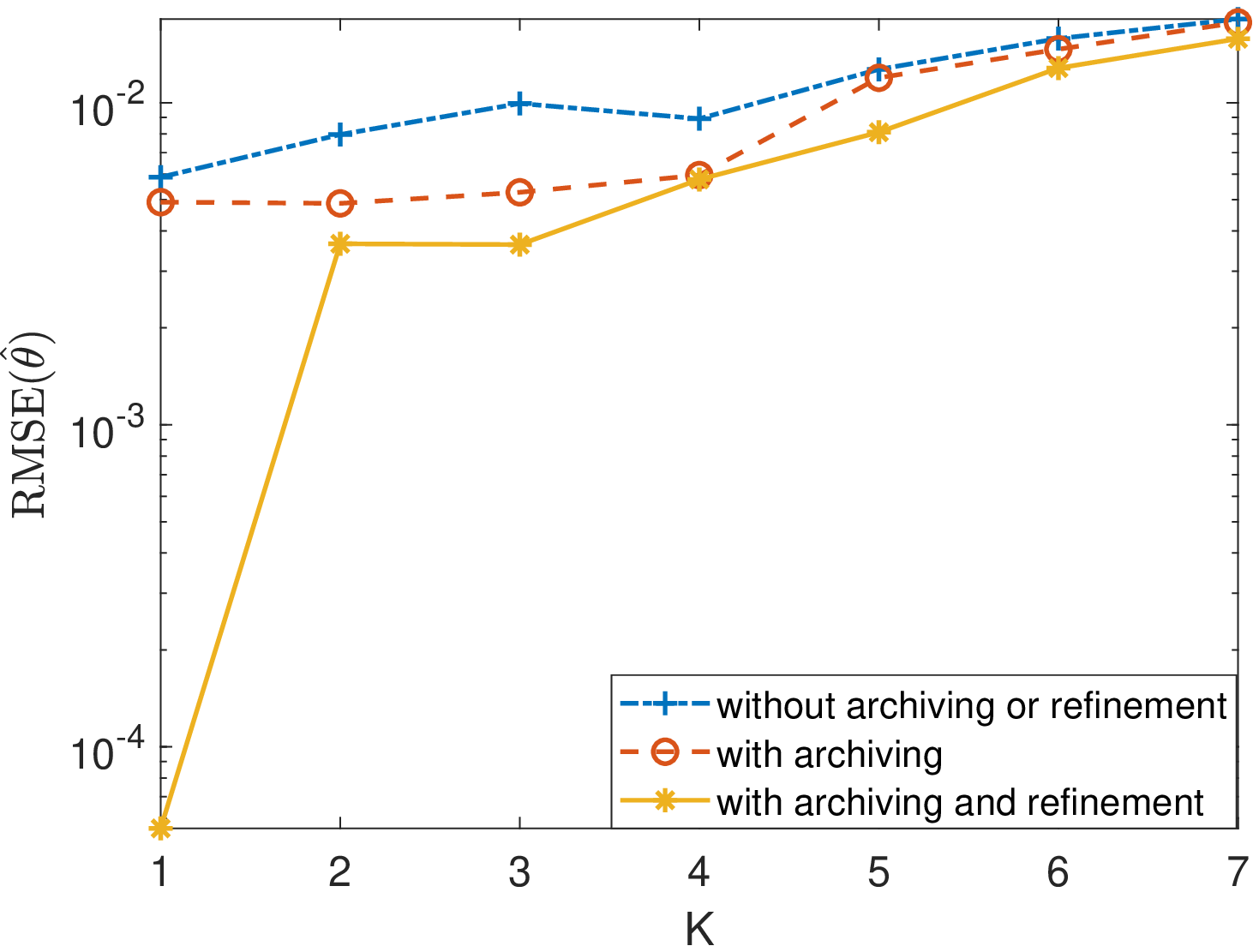}} \qquad
	\subfigure[]{\includegraphics[width=0.42\textwidth]{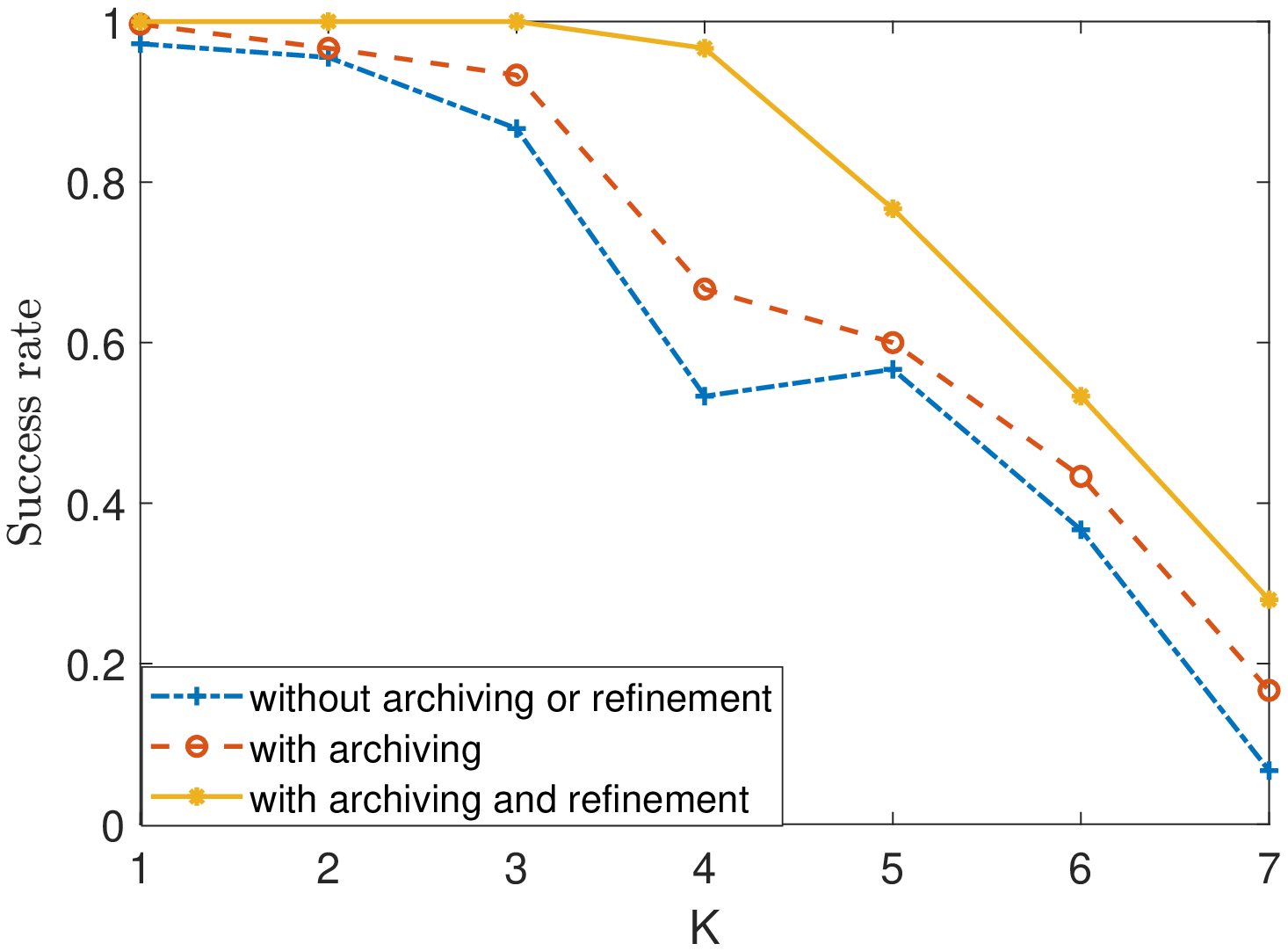}}
	\caption{RMSE($\hat{\mta}$) (a) and success rate (b) results of three versions of MVESA versus model order $K$ for $M=15$, $T=10$, and SNR$=10$dB.} 
	\label{fig-operators}
\end{figure*}

\begin{figure*}[t] 
	\centering
	\subfigure[]{\includegraphics[width=0.42\textwidth]{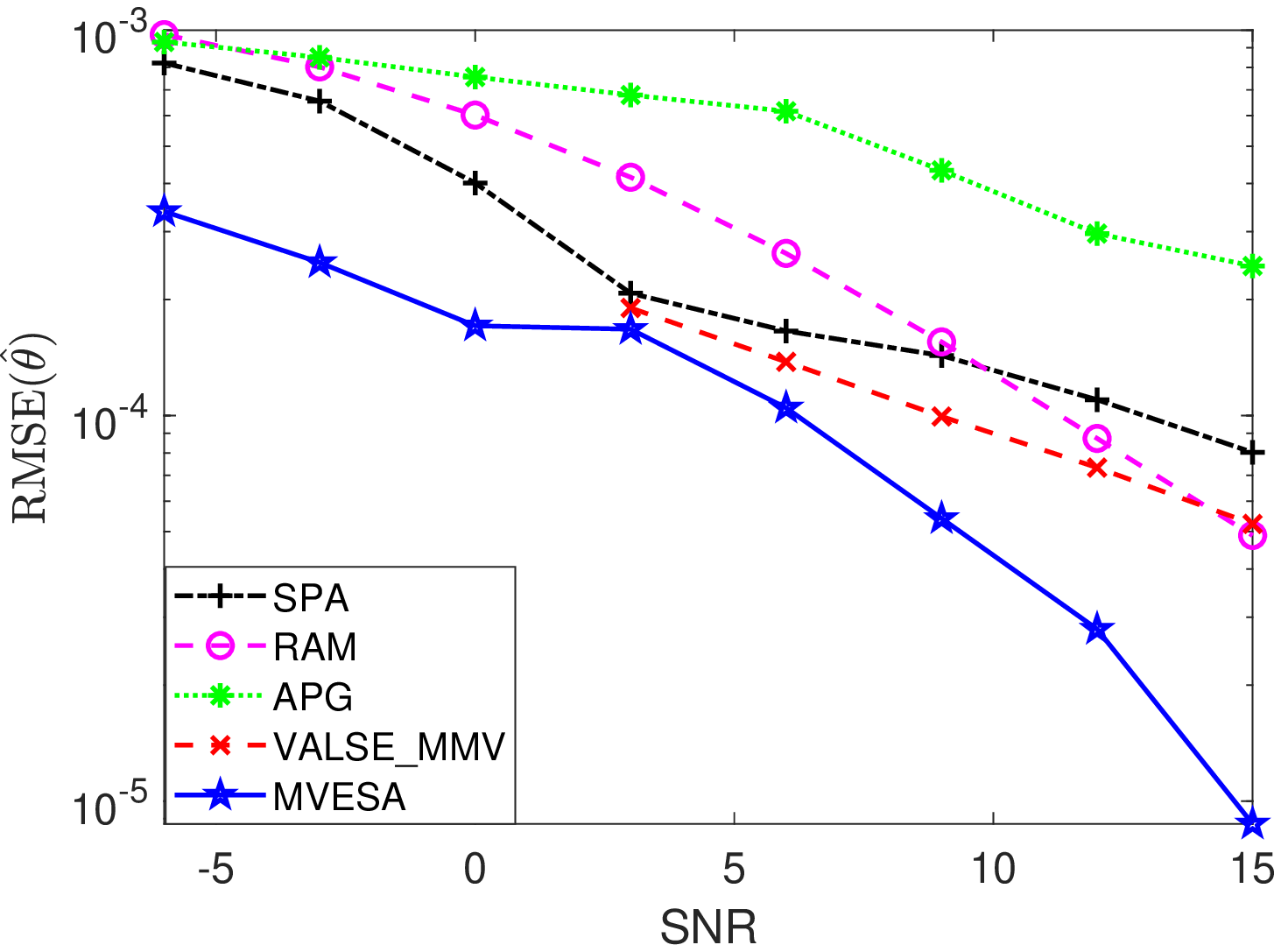}} \qquad
	\subfigure[]{\includegraphics[width=0.42\textwidth]{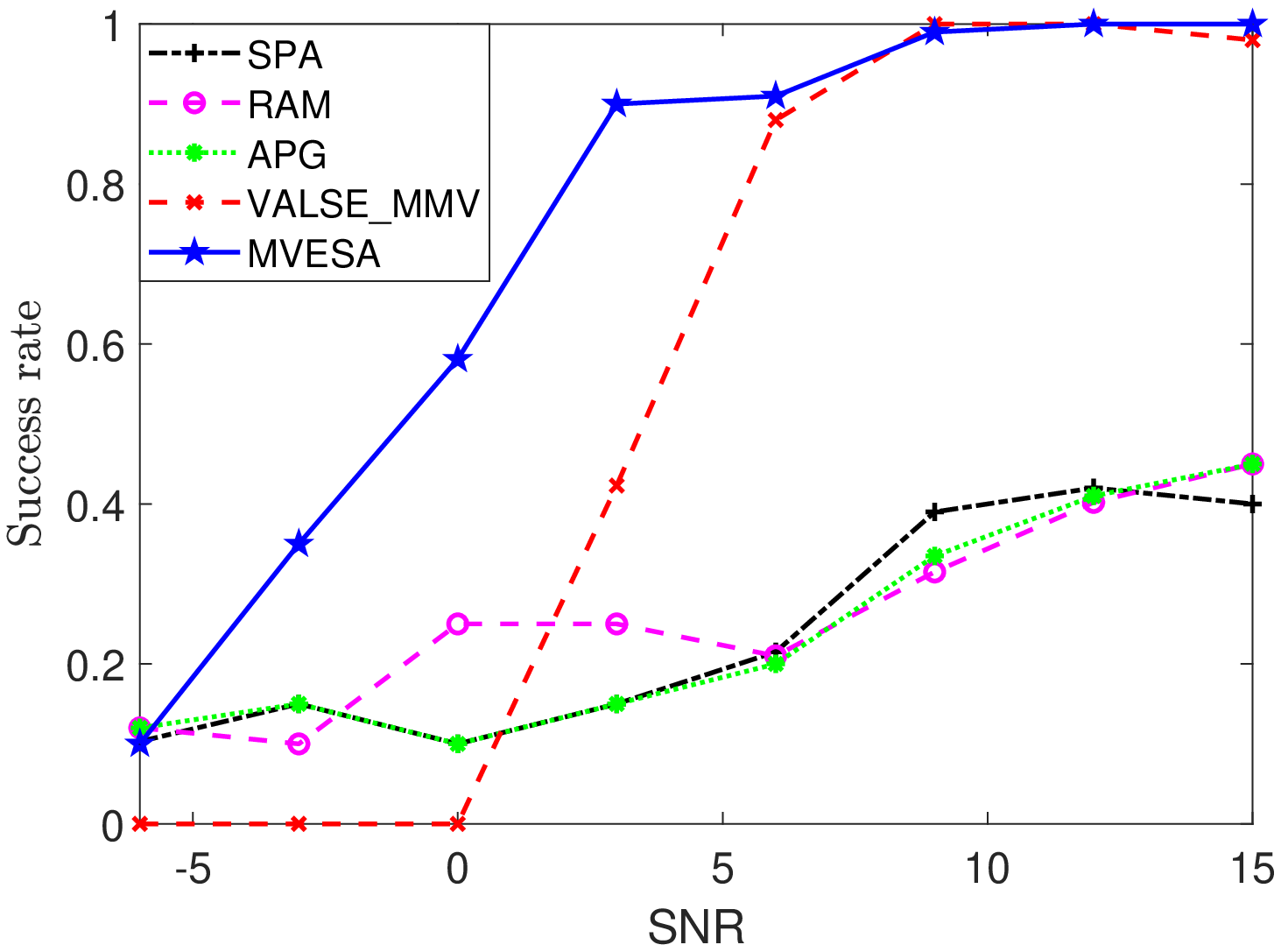}}
	\caption{RMSE($\hat{\mta}$) (a) and success rate (b) results of all algorithms versus SNR for $K=4$, $M=15$, and $T=30$.} 
	\label{fig-SNR}
\end{figure*}

\begin{figure*}[] 
	\centering
	\subfigure[]{\includegraphics[width=0.42\textwidth]{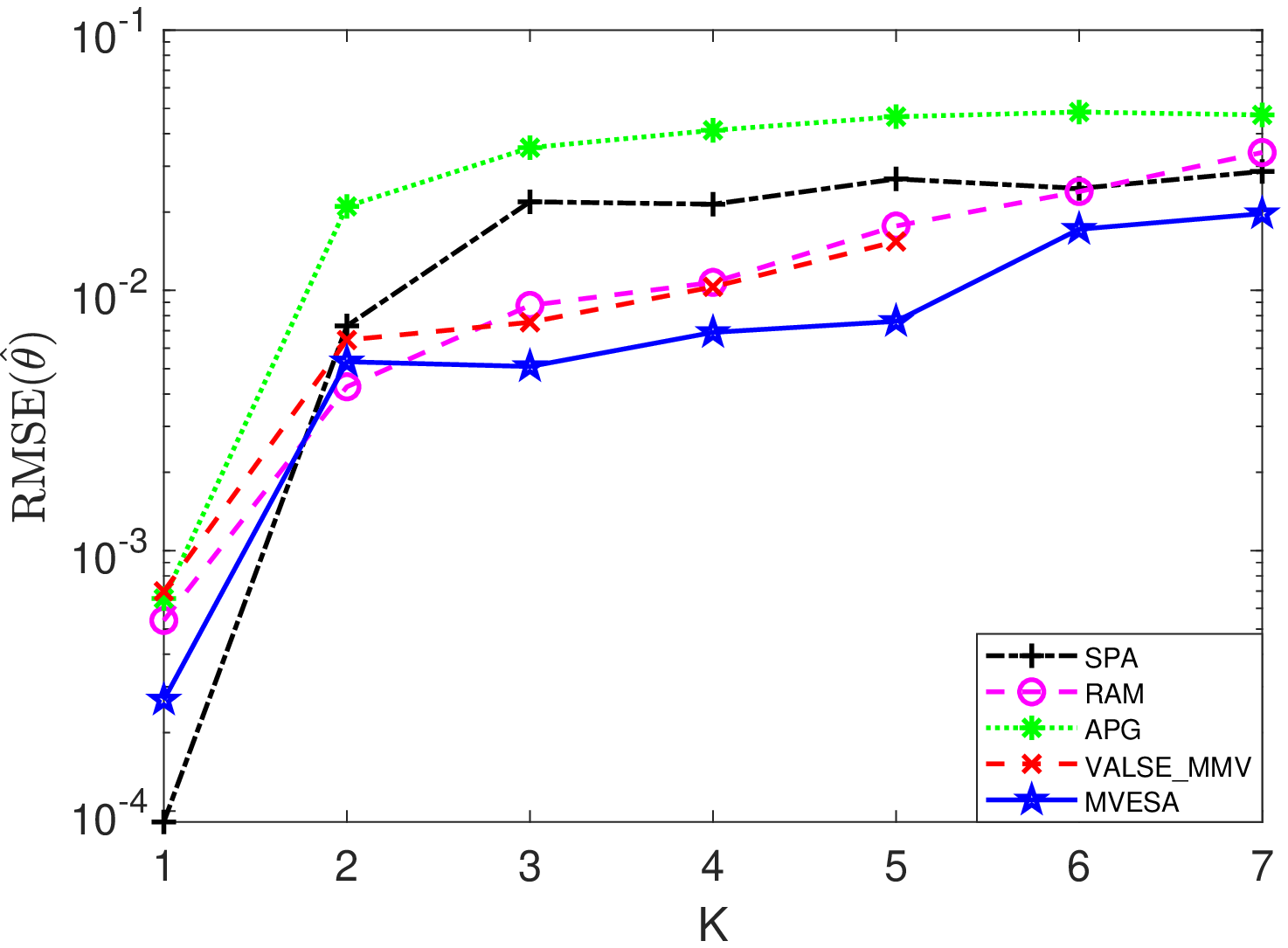}} \qquad
	\subfigure[]{\includegraphics[width=0.42\textwidth]{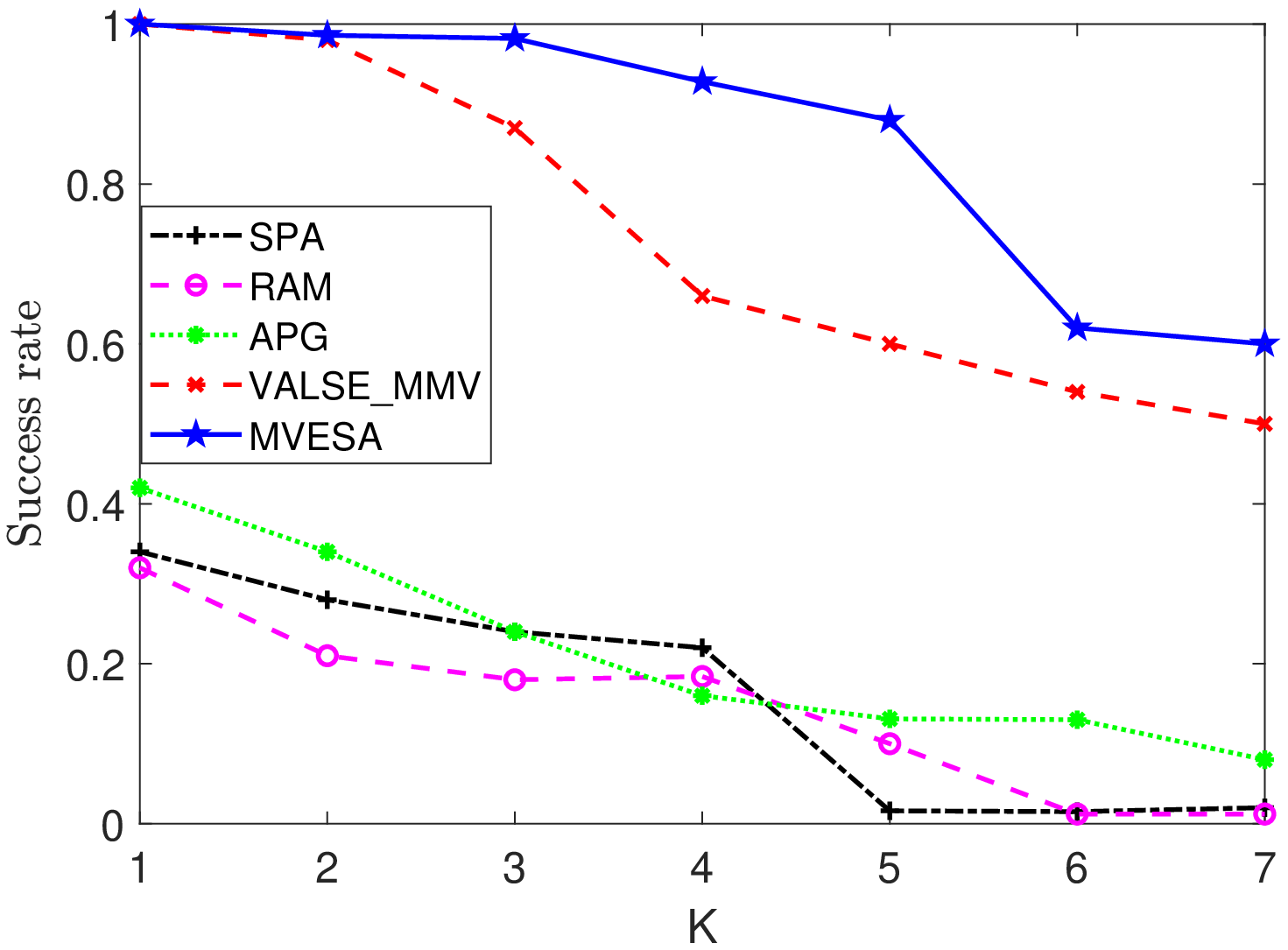}}
	\caption{RMSE($\hat{\mta}$) (a) and success rate (b) results of all algorithms versus model order $K$ for $M=15$, $T=10$ and SNR$=10$dB.} %thetas=[-22.6  -19.7  -13.4 19.9 36.6  45.1 56.8];
	\label{fig-K}
\end{figure*}

\subsection{Setup, Metrics and Algorithms}
\textbf{Setup}. According to the LSE model (\ref{eq-y=As+noise}), a number of $K$ frequencies are randomly generated within $[-1,1)$. The amplitudes $\St$ are drawn from i.i.d. from $\mathcal{CN}(1, 0.1)$. Note that we do not control the minimum frequency separation, thus the frequencies may not be guaranteed to be recovered, even for large $M$. The noise samples contaminating the measurements are independent and zero-mean complex Gaussian distributed. 

\textbf{Metrics}. Since Bayesian-based methods and our algorithm do not output spatial spectrum, spatial spectrum will not be used for comparison. We employ two statistical measures, i.e., root mean square error ($\text{$\text{RMSE}$}$) and success rate. RMSE is obtained by averaging the frequency combination error over $\upsilon$ Monte Carlo runs:
\begin{equation} 
\label{eq-RMSE}
\begin{aligned}
\text{RMSE}(\hat{\mta})=\sqrt{\frac{1}{\upsilon}\sum_{i=1}^{\upsilon}\|\hat{\mta}-\mta\|_2},
\end{aligned}
\end{equation} 
where $\hat{\mta}$ and $\mta$ are the estimated and true frequency combination, respectively. The averaging operation performs over the trials in which the estimated model number is greater than or equal to $K$. The assignment of estimated frequencies to the true one is executed based on the Hungarian algorithm \cite{1962Algorithms}. The success rate is the empirical probability that the estimated model order $\hat{K}$ is the same as the true value $K$, denoted as $Pr(K=\hat{K})$. % For a given simulation point, the RMSE is obtained by averaging only the trials in which the estimated source number is greater than or equals to $K$.

\textbf{Algorithms}. We conduct simulations to compare the performance of MVESA with the state-of-the-art gridless algorithms, i.e., SPA \cite{yang2014a}, RAM \cite{yang2015enhancing}, APG \cite{wagner2019gridless} and VALSE\_MMV \cite{zhu2019grid}. The comparison algorithms are introduced below. 
\begin{itemize}
	\item SPA: A gridless algorithm based on covariance fitting criteria and convex optimization. This algorithm can work without model order but cannot determine it accurately.
	\item RAM: A gridless algorithm based on reweighted atomic-norm minimization for enhancing sparsity and resolution. The model order is exploited by reweighted atomic norm.	
	\item APG: A gridless algorithm that directly solves the atomic $l_0$ norm minimization problem via alternating projections. But the model order needs to be known as a priori. 
	\item VALSE\_MMV: A representative gridless sparse Bayesian inference-based algorithm that estimates the posterior probability density functions of frequencies. The model order is estimated using Bernoulli-Gaussian distribution.	
\end{itemize} 

Comparison algorithms' parameters are set in accordance with their original papers  \cite{yang2014a}\cite{yang2015enhancing}\cite{wagner2019gridless}\cite{zhu2019grid}, respectively. For SPA, RAM and APG, the model order is set to its possible maximum value, $M-1$. For proposed MVESA, we set population size$=30$, mutation distribution index$=20$, and mutation probability$=1/k$, where $k$ is the model order of current solution. To accelerate the search efficiency of MVESA, we generate the initial population in this way: a single solution with a maximum length $M-1$ is produced by the simple Capon method \cite{2005MalioutovA}, and the remaining $N-1$ solutions are randomly geneated.

For a fair comparison, all the algorithms stop running when the change of estimated measurements $\|\hat{\Yt}^G-\hat{\Yt}^{G-1}\|_F/\|\hat{\Yt}^{G-1}\|_F$ is less than $10^{-6}$ in three consecutive generations, or the total number of iterations for comparison methods and our method exceeds 5000 and 100 respectively. The reason of setting the number of iterations like this is that, at each generation, no more than 50 solutions are explored in our method, and one solution explored in comparison methods. As a result, all algorithms are allowed to explore at most 5000 solutions in total, providing fair comparison. The total number of Monte Carlo runs are set to 200 of for all algorithms. 

\subsection{Detailed Analysis of MVESA} \label{sec-exp-BEA}
In this subsection, the effectiveness of the two-objective function and the proposed archiving and model order pruning mechanism are investigated to demonstrate the superiority of MVESA.

\subsubsection{Study of Objective function} \label{sec-exp-knee}
The two proposed objectives, including the model order and measurement error, are conflicting with each other. This conflicting characteristic enables MVESA to determine the model order automatically. To validate it, we conduct a simulation with $K=4$, $M=15$, $T=20$ and increase SNR from -5dB to 15 dB. Apart from this, the noiseless case is involved as reference. Fig. \ref{fig-PF-and-slope} depicts the typical Pareto front (a) and slope variance (b) results of the final archive over 200 runs. Fig. \ref{fig-PF-and-slope}(b) is obtained by computing the slope variance of the Pareto front according to the kink method \cite{mierswa2006information}. The knee solutions characterized by the maximum slope variance for different SNRs are identified and enframed within dotted line in Fig. \ref{fig-PF-and-slope}(a). It can be seen that, the identified knee solutions provide the best trade-off between the two conflicting objectives and acquire the true model order. Thus, it indicates the efficiency of the two-objective function of MVESA.

\begin{figure*}[] 
	\centering
	\subfigure[]{\includegraphics[width=0.42\textwidth]{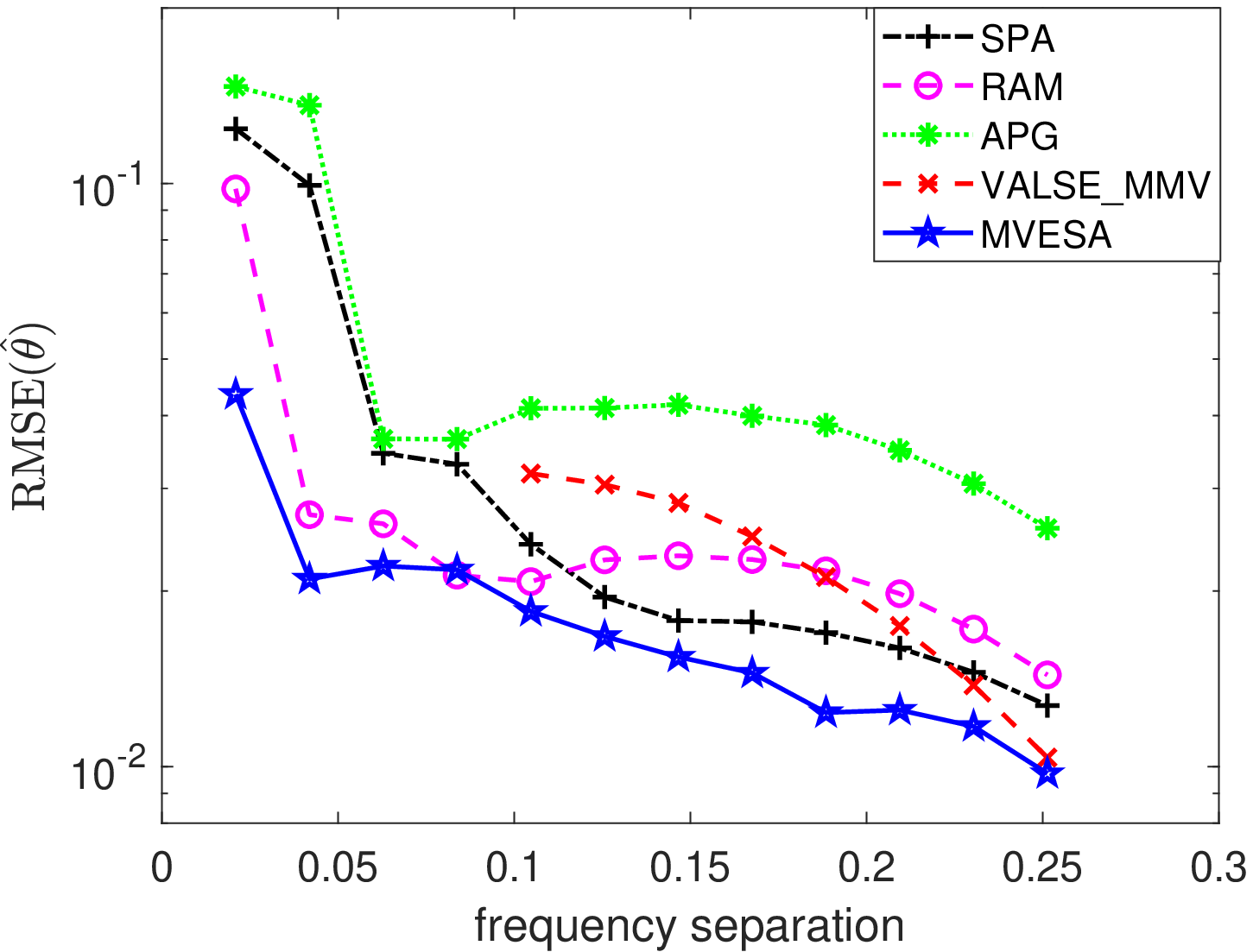}} \qquad
	\subfigure[]{\includegraphics[width=0.42\textwidth]{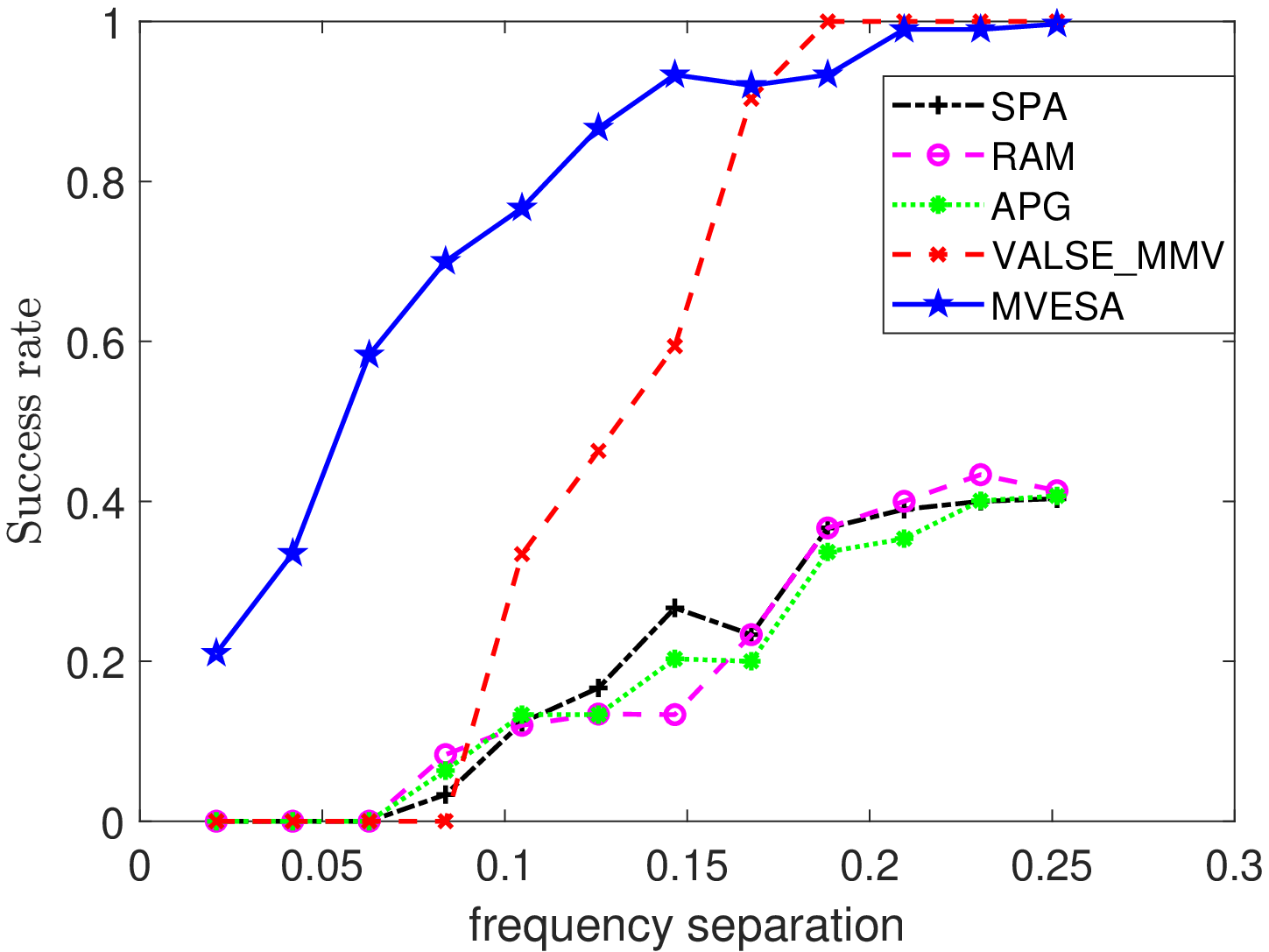}}
	\caption{RMSE($\hat{\mta}$) (a) and success rate (b) results of all algorithms versus frequency separation for $M=6$, $T=10$, SNR$=10$dB.} %thetas=[-22.6  -19.7  -13.4 19.9 36.6  45.1 56.8];
	\label{fig-separation}
\end{figure*}

\begin{figure*}[] 
	\centering
	\subfigure[]{\includegraphics[width=0.42\textwidth]{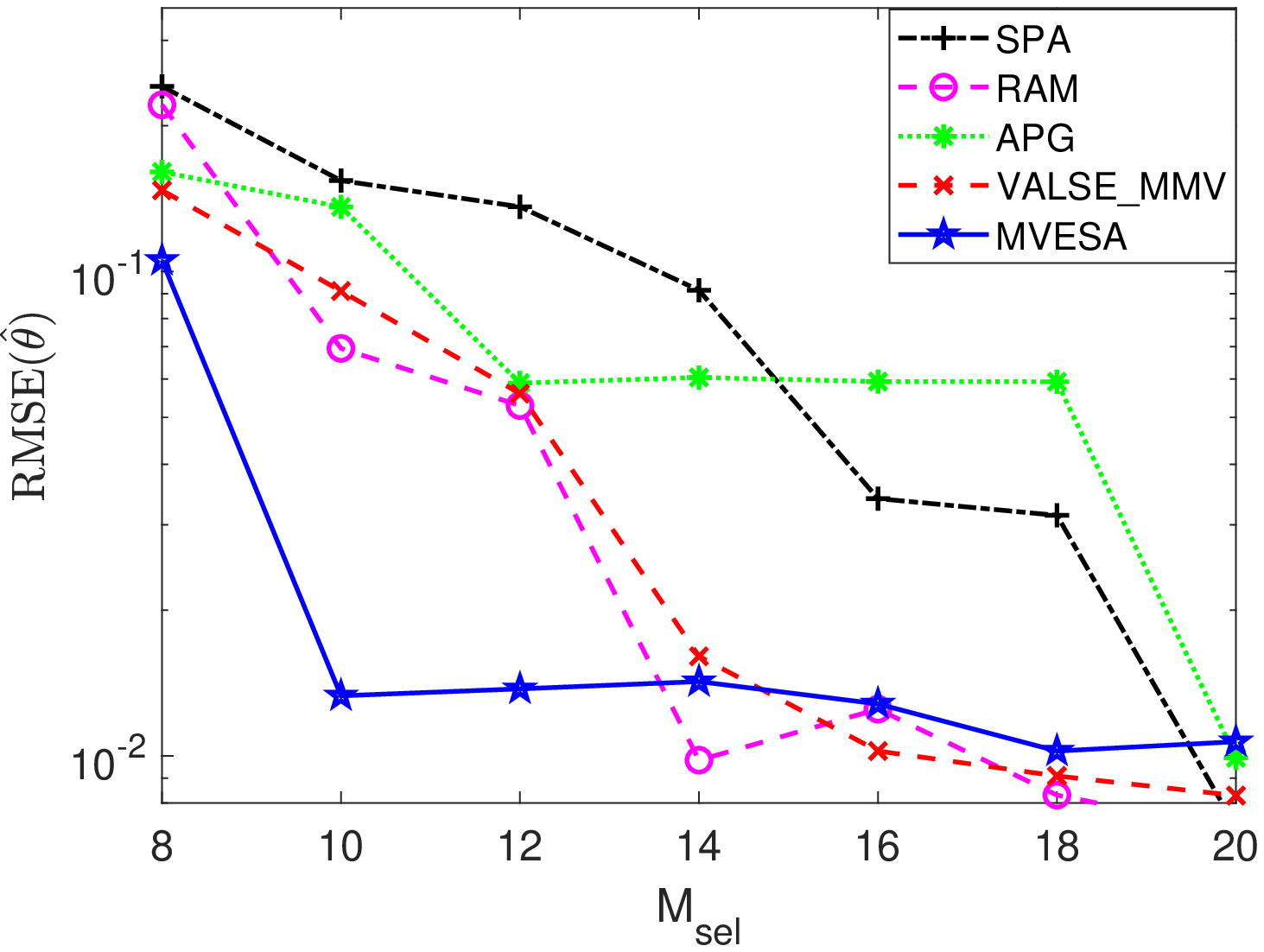}} \qquad
	\subfigure[]{\includegraphics[width=0.42\textwidth]{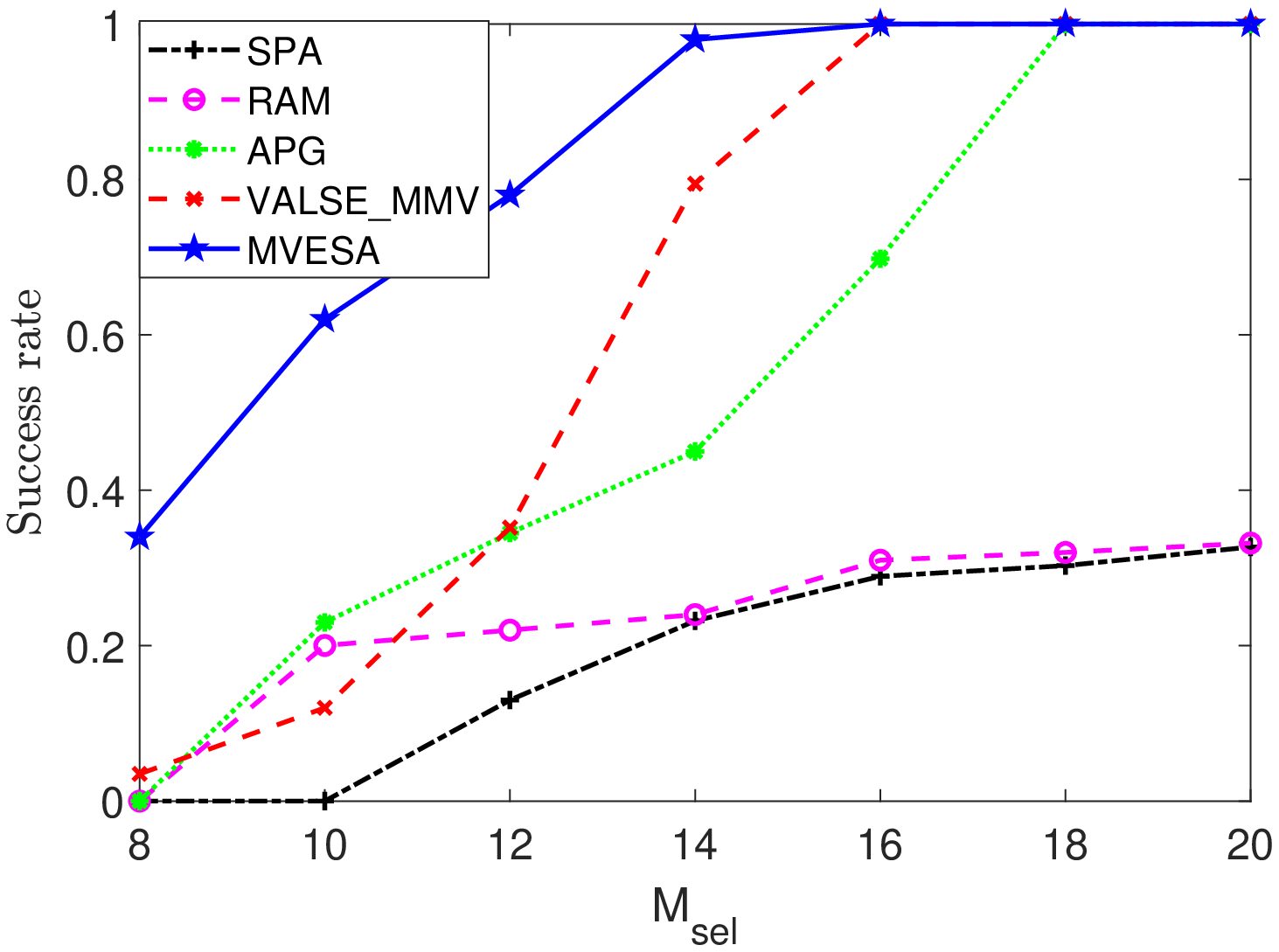}}
	\caption{RMSE($\hat{\mta}$) (a) and success rate (b) results of all algorithms under varying $M_{sel}$s for $K=3$, $T=10$ and SNR$=10$dB.} %thetas=[-22.6  -19.7  -13.4 19.9 36.6  45.1 56.8];
	\label{fig-M}
\end{figure*}

\subsubsection{Study of the archiving and model order pruning mechanism} 
To validate this mechanism's efficacy, we compare the performance of three versions of MVESA, including the first version without archiving or model order pruning, the second version with archiving, and the third version with archiving and model order pruning. The simulation parameters are set to $M=15$, $T=10$ snapshots, SNR$=10$dB, and the number frequencies $K$ increases from 1 to 7. Fig. \ref{fig-operators} plots the resultant RMSEs and success rates under different $K$s. The third version, i.e., MVESA, is always superior to the other versions in terms of frequency error and success rate. This advantage can be explained as follows. Compared to the first version, the last two ones incorporate archive to collect the best solutions so far. Consequently, they can avoid missing the optimal solution, providing better estimation performance. The performance gap between the second version and MVESA is because the model order pruning operation can fully explore the search space over different dimensionality. This operation could skip suboptimal solutions, bring enhanced convergence performance. 

\subsection{Comparison of BEA Against Other Methods} \label{sec-exp-comparison}
In this subsection, the algorithms' ability of handling complete data and incomplete data are investigated, respectively. The complete data is referred to as all $M$ measurements being available while incomple data mean that some of the $M$ measurements are missing. The missing data case may occur due to sensor failure, outliers, or other physical constraints. The time cost of all these algorithms are also compared.

\subsubsection{Handling complete data}
In Simulation 1, Monte Carlo trials are carried out to investigate the robustness to noise intensity. We set the parameters $K=4$, $M=15$, $T=30$ snapshots, and vary SNRs from -6dB to 15dB. Fig. \ref{fig-SNR} shows the RMSE and success rate and frequency errors versus SNRs. MVESA retains much lower frequency errors compared to other algorithms for all SNRs. In terms of success rate, SPA, RAM and APG roughly fail to determine the model order. By contrast, VALSE\_MMV and MVESA estimate the model order more accurately. Note that MVESA obtains the highest success rate at very low SNRs. The superior performance of MVESA is because it owns outstanding ability of exploring the dimensionality-changing search space, which can skip suboptimal solutions and bring enhanced performance.    

In Simulation 2, the capability of model order selection in scenarios with varying model orders is investigated. The parameters are set to $M=15$, $T=10$ snapshots, SNR$=10$dB, and model order $K$ increasing from 1 to 7. Results of RMSE and success rate under varying model orders are shown in Fig. \ref{fig-K}. SPA, RAM and APG still perform worse than VALSE\_MMV and MVESA both in terms of RMSE and success rate, because they lack the capability of model order selection. MVESA shows an absolute advantage over VALSE\_MMV in terms of the two evaluation metrics for most of $K$s. This advantage benefits from the atomic $l_0$ norm, which helps predict the model order more accurately and largely improve the estimation accuracy.

Simulation 3 studies the statistical performance of LSE of resolving two closely-spaced frequencies. Assume the distance between two components vary from 0.02 to 0.26, $M=6$, $T=10$ snapshots, and SNR$=10$dB, Fig. \ref{fig-separation} shows the results versus frequency separation. It can be observed that MVESA wins the best in 10 out of the 12 cases in terms of frequency error. When the two frequencies are located relatively closely (separation between 0.04 and 0.15), MVESA obtains a significantly high resolution, while other algorithms almost fail to work due to the resolution limit caused by suboptimal sparse metrics. The high resolution of MVESA validates the importance of the atomic $l_0$ norm for promoting sparsity.

\subsubsection{Handling on incomplete data}
Simulation 4 investigates the capability of handling incomplete sample data. We consider the estimation of $M=20$, $K=3$, $T=10$ snapshots, and SNR$=10$dB. Particularly, we extract $M_{sel}\leq M$ measurements from complete measurements $\Yt$ with indices in $\mathcal{M}\subseteq\{0,...M-1\}$, $|\mathcal{M}|=M_{sel}$, thus the resultant measurments data are incomplete. The RMSE and success rate results for incomplete data are shown in Fig. \ref{fig-M}. It can be observed that for $M_{sel}<14$, MVESA retains a better estimation performance in terms of RMSE and success rate. For $14\leq M_{sel}<20$, MVESA achieves slightly worse frequency error than VALSE\_MMV and RAM, and far surpass SPA and APG. This phenonmenon may be because MVESA ignores the noise in modeling and is not statistically inconsistent in $M$. Regardless of this, MVESA always achieves the right model order with highest probability for all $M_{sel}s$, which validates the effective of MVESA in joint estimation of frequencies and model order.   

\subsubsection{Time cost analysis}
Simulate 5 compares the time complexity of all algorithms under different $M$s. This simulation is implemented in MATLAB R2018b on a PC with Intel i7-7700 CPU and 32GB RAM. Parametric settings are set the same as Fig. \ref{fig-M}. The computational time is displayed in Table \ref{Tab-time}. VALSE\_MMV runs the fastest, following by SPA. MVESA is slightly slower than SPA, but outdistances RAM and APG. This is because MVESA needs to execute matrix inversions, but it does not require to solve semidefinite programming problems. Considering the parallel nature of evolutionary algorithms \cite{2015Distributed}, we suggest to accelerate MVESA by parallel implementation to satisfy large-scaled real-world applications.  

\begin{table}[t]
	\renewcommand{\arraystretch}{1.45}
	\centering
	\caption{Average running time (in seconds) of all algorithms versus $M$}
	\label{Tab-time}
	\begin{tabular}{ccccccc}
		\hline
		$M$ & SPA & RAM  & APG &  VALSE\_MMV  & MVESA \\ \hline
		8   & 0.7609   & 5.3165   & 0.8966   & 0.0400   & 0.6170 \\
		12  & 0.7846   & 5.6782   & 1.5747   & 0.0603   & 0.7068 \\
		16  & 0.8237   & 5.7789   & 2.6863   & 0.0913   & 0.8037 \\
		20  & 0.9202   & 5.6872   & 4.5730   & 0.1133   & 0.9575 \\
		24  & 0.9233   & 5.8040   & 6.5828   & 0.1443   & 1.0924 \\ \hline
	\end{tabular}%} % 
\end{table}  

\section{Conclusion}
In this paper, we have proposed a novel idea of simultaneously estimating the model order and frequencies by means of atomic $l_0$ norm. To accomplish this, we have built a multiobjective optimization model, with the measurement error and atomic $l_0$ norm being the two objectives. The atomic $l_0$ norm directly exploits sparsity without relaxations, breaking the resolution limit and estimating the model order accurately. To solve the resultant NP-hard problem, we have designed the multiobjective variable-length evolutionary search algorithm with two innovations. One is the variable-length coding and search strategy, which provides a flexible representation of frequencies with different sizes, and implement full exploration over the variable and open-ended search space. Another innovation is the model order pruning mechanism, which reduces the solutions' redundancy by heuristically pruning less contributive frequencies. This mechanism highly improves the convergence and diversity performance. Experiments results have demonstrated the superiority of MVESA in terms of RMSE and success rate.

The proposed method involves matrix inversions, which is time-consuming for large-scaled LSE problems. Therefore, in the future, we plan to design more computational efficient methods. We also expect to further improve MVESA to be statistically consistent so that the frequency estimation performance in cases with large $M$s would be improved.

%\section*{Acknowledgment}

\ifCLASSOPTIONcaptionsoff
  \newpage
\fi


\bibliographystyle{IEEEtran}
\bibliography{GridlessRef}

\iffalse
\begin{IEEEbiography}{Michael Shell}
Biography text here.
\end{IEEEbiography}

\begin{IEEEbiographynophoto}{John Doe}
Biography text here.
\end{IEEEbiographynophoto}

%\newpage

\begin{IEEEbiographynophoto}{Jane Doe}
Biography text here.
\end{IEEEbiographynophoto}
\fi

%\vfill


\end{document}